\documentclass[12pt]{article}
\usepackage{amssymb}
\usepackage{graphicx}
\usepackage{latexsym}
\usepackage{tikz}

\usepackage{amsmath}

\def\part#1{\frac{\partial\phantom{q}}{\partial#1}}

\newenvironment{rmk}{\begin{trivlist}\item[]{\bf Remark:} }
{\end{trivlist}}
\newenvironment{rmks}{\begin{trivlist}\item[]{\bf Remarks:} }
{\end{trivlist}}
\newenvironment{ex}{\begin{trivlist}\item[]{\bf Example:} }
{\end{trivlist}}

\newenvironment{prf}{\begin{trivlist}\item[]{\bf Proof:} }
{\hfill $\Box$ \end{trivlist}}

\newtheorem{thm}{Theorem}

\newtheorem{prp}[thm]{Proposition}

\newcommand{\lie}[1]{\mathfrak{#1}}
\def\End{\mathop{\rm End}\nolimits}

\font\dynkfont=cmsy10 scaled\magstep5    \skewchar\dynkfont='60
\def\dynk{\textfont2=\dynkfont}
\def\hr#1,#2;{\dimen0=.4pt\advance\dimen0by-#2pt
              \vrule width#1pt height#2pt depth\dimen0}
\def\vr#1,#2;{\vrule height#1pt depth#2pt}
\def\blb#1#2#3#4#5
            {\hbox{\ifnum#2=0\hskip11.5pt
                   \else\ifnum#2=1\hr13,5;\hskip-1.5pt
                   \else\ifnum#2=2\hr13.5,6.5;\hskip-13.5pt
                                  \hr13.5,3.5;\hskip-2pt
                   \else\ifnum#2=3\hr13.7,8;\hskip-13.7pt
                                  \hr13,5;\hskip-13pt
                                  \hr13.7,2;\hskip-2.2pt\fi\fi\fi\fi
                   $#1$
                   \ifnum#4=0
                   \else\ifnum#4=1\hskip-9.2pt\vr22,-9;\hskip8.8pt
                   \else\ifnum#4=2\hskip-10.9pt\vr22,-8.75;\hskip3pt
                                  \vr22,-8.75;\hskip7.1pt
                   \else\ifnum#4=3\hskip-12.6pt\vr22,-8.5;\hskip3pt
                                  \vr22,-9;\hskip3pt
                                  \vr22,-8.5;\hskip5.4pt\fi\fi\fi\fi
                   \ifnum#5=0
                   \else\ifnum#5=1\hskip-9.2pt\vr1,12;\hskip8.8pt
                   \else\ifnum#5=2\hskip-10.9pt\vr1.25,12;\hskip3pt
                                  \vr1.25,12;\hskip7.1pt
                   \else\ifnum#5=3\hskip-12.6pt\vr1.5,12;\hskip3pt
                                  \vr1,12;\hskip3pt
                                  \vr1.5,12;\hskip5.4pt\fi\fi\fi\fi
                   \ifnum#3=0\hskip8pt
                   \else\ifnum#3=1\hskip-5pt\hr13,5;
                   \else\ifnum#3=2\hskip-5.5pt\hr13.5,6.5;
                                  \hskip-13.5pt\hr13.5,3.5;
                   \else\ifnum#3=3\hskip-5.7pt\hr13.7,8;
                                  \hskip-13pt\hr13,5;
                                  \hskip-13.7pt\hr13.7,2;\fi\fi\fi\fi
 
                   }}
\def\blob#1#2#3#4#5#6#7{\hbox
{$\displaystyle\mathop{\blb#1#2#3#4#5 }_{#6}\sp{#7}$}}
\def\up#1#2{\dimen1=33pt\multiply\dimen1by#1\hbox{\raise\dimen1\rlap{#2}}}
\def\uph#1#2{\dimen1=17.5pt\multiply\dimen1by#1\hbox{\raise\dimen1\rlap{#2}}}
\def\dn#1#2{\dimen1=33pt\multiply\dimen1by#1\hbox{\lower\dimen1\rlap{#2}}}
\def\dnh#1#2{\dimen1=17.5pt\multiply\dimen1by#1\hbox{\lower\dimen1\rlap{#2}}}
 
\def\rlbl#1{\kern-8pt\raise3pt\hbox{$\scriptstyle #1$}}
\def\llbl#1{\raise3pt\llap{\hbox{$\scriptstyle #1$\kern-8pt}}}
\def\elbl#1{\kern3pt\lower4.5pt\hbox{$\scriptstyle #1$}}
\def\lelbl#1{\rlap{\hbox{\kern-9pt\raise2.5pt\hbox{{$\scriptstyle #1$}}}}}

\def\whtd#1#2#3#4#5{\blob\circ#1#2#3#4{#5}{}}
\def\blkd#1#2#3#4#5{\blob\bullet#1#2#3#4{#5}{}}
\def\whtu#1#2#3#4#5{\blob\circ#1#2#3#4{}{#5}}
\def\blku#1#2#3#4#5{\blob\bullet#1#2#3#4{}{#5}}
\def\whtr#1#2#3#4#5{\blob\circ#1#2#3#4{}{}\rlbl{#5}}
\def\blkr#1#2#3#4#5{\blob\bullet#1#2#3#4{}{}\rlbl{#5}}

\def\rwng{\hbox{$\vbox{\offinterlineskip{
  \hbox{\phantom{}\kern6pt{$\circ$}}\kern-2.5pt\hbox{$\Biggr/$}\kern-0.5pt
  \hbox{\phantom{}\kern-5pt$\circ$}\kern-3.0pt\hbox{$\Biggr\backslash$}
  \kern-1.5pt\hbox{\phantom{}\kern6pt{$\circ$}} }}$}}
 
\def\lwng{\hbox{$\vbox{\offinterlineskip{ \hbox{$\circ$}
  \kern-3.0pt\hbox{\phantom{}\kern6.0pt{$\Biggr\backslash$}}
  \kern-0.5pt\hbox{\phantom{}\kern11pt{$\circ$}}\kern-3.5pt
  \hbox{\phantom{}\kern5.0pt {$\Biggr/$}}\kern-1.0pt\hbox{$\circ$} }}$}}

\def\drwng#1#2#3{\hbox{$\vcenter{ \offinterlineskip{
  \hbox{\phantom{}\kern7pt{$\circ^{\elbl{#3}}$}}
  \kern-2.5pt\hbox{$\Biggr/$}\kern-0.5pt
  \hbox{\phantom{}\kern-5pt$\circ^{ \elbl{#1}}$}
  \kern-3.0pt\hbox{$\Biggr\backslash$}
  \kern-1.5pt\hbox{\phantom{}\kern7pt{$\circ^{\elbl{#2}}$}}  } }$}}
\def\drwngt#1#2#3{\hbox{$\vcenter{ \offinterlineskip{
  \hbox{\phantom{}\kern7pt{$\bullet^{\elbl{#3}}$}}
  \kern-2.5pt\hbox{$\Biggr/$}\kern-0.5pt
  \hbox{\phantom{}\kern-5pt$\circ^{ \elbl{#1}}$}
  \kern-3.0pt\hbox{$\Biggr\backslash$}
  \kern-1.5pt\hbox{\phantom{}\kern7pt{$\bullet^{\elbl{#2}}$}}  } }$}}
 
\def\dlwng#1#2#3{\hbox{$\vcenter{\offinterlineskip{ \hbox{$\lelbl{#1}\circ$}
  \kern-3.0pt\hbox{\phantom{}\kern6.0pt{$\Biggr\backslash$}}
  \kern-0.5pt\hbox{\phantom{}\kern11pt{$\lelbl{#2}\circ$}}\kern-3.5pt
  \hbox{\phantom{}\kern5.0pt {$\Biggr/$}}\kern-1.0pt\hbox{$\lelbl{#3}\circ$}}}$}
 }

\def\rde#1#2#3{\hbox{\phantom{}\kern-4pt\hbox{$\vcenter{\offinterlineskip  \hbox
{
               \raise 4.5pt\hbox{\vrule height0.4pt width13pt depth0pt}
                \kern-1pt\vbox{ \hbox{\drwng{#1}{#2}{#3}}} }}$  }}  }
\def\rdet#1#2#3{\hbox{\phantom{}\kern-4pt\hbox{$\vcenter{\offinterlineskip \hbox
{
               \raise 4.5pt\hbox{\vrule height0.4pt width13pt depth0pt}
                \kern-1pt\vbox{ \hbox{\drwngt{#1}{#2}{#3}}} }}$  }}  }
 
\def\lde#1#2#3{\hbox{$\vcenter{\offinterlineskip  \hbox{
               \dlwng{#1}{#2}{#3}\kern-4.2pt\lower0.4pt\hbox{$\vcenter{\hrule 
                               width13pt}$}
               \kern-8pt\phantom{}   }}  $}}

\def\rwngb{\hbox{$\vbox{\offinterlineskip{
  \hbox{\phantom{}\kern6pt{$\bullet$}}\kern-2.5pt\hbox{$\Biggr/$}\kern-0.5pt
  \hbox{\phantom{}\kern-5pt$\bullet$}\kern-3.0pt\hbox{$\Biggr\backslash$}
  \kern-1.5pt\hbox{\phantom{}\kern6pt{$\bullet$}} }}$}}
 
\def\lwngb{\hbox{$\vbox{\offinterlineskip{ \hbox{$\bullet$}
  \kern-3.0pt\hbox{\phantom{}\kern6.0pt{$\Biggr\backslash$}}
  \kern-0.5pt\hbox{\phantom{}\kern11pt{$\bullet$}}\kern-3.5pt
  \hbox{\phantom{}\kern5.0pt {$\Biggr/$}}\kern-1.0pt\hbox{$\bullet$} }}$}}

\def\dbrwng#1#2#3{\hbox{$\vcenter{ \offinterlineskip{
  \hbox{\phantom{}\kern6pt{$\bullet^{\elbl{#3}}$}}
  \kern-2.5pt\hbox{$\Biggr/$}\kern-0.5pt
  \hbox{\phantom{}\kern-5pt$\bullet^{ \elbl{#1}}$}
  \kern-3.0pt\hbox{$\Biggr\backslash$}
  \kern-1.5pt\hbox{\phantom{}\kern6pt{$\bullet^{\elbl{#2}}$}}  } }$}}
 
\def\dblwng#1#2#3{\hbox{$\vcenter{\offinterlineskip{ \hbox{$\lelbl{#1}\bullet$}
  \kern-3.0pt\hbox{\phantom{}\kern6.0pt{$\Biggr\backslash$}}
  \kern-0.5pt\hbox{\phantom{}\kern11pt{$\lelbl{#2}\bullet$}}\kern-3.5pt
  \hbox{\phantom{}\kern5.0pt {$\Biggr/$}}\kern-1.0pt\hbox{$\lelbl{#3}\bullet$}}}
$} }

\def\rbde#1#2#3{\hbox{\phantom{}\kern-4pt\hbox{$\vcenter{\offinterlineskip  \hbo
x{
               \raise 4.5pt\hbox{\vrule height0.4pt width13pt depth0pt}
                \kern-1pt\vbox{ \hbox{\dbrwng{#1}{#2}{#3}}} }}$  }}  }
 
\def\lbde#1#2#3{\hbox{$\vcenter{\offinterlineskip  \hbox{
               \dblwng{#1}{#2}{#3}\kern-4.2pt\lower0.4pt\hbox{$\vcenter{\hrule w
idth13pt}$}
               \kern-8pt\phantom{}   }}  $}}

\def\ddgu#1.#2.{\dynk  \whtu0300{#1}\blku3000{#2}}
\def\ddgd#1.#2.{\dynk  \whtd0300{#1}\blkd3000{#2}}
 
\def\eddgiu#1.#2.#3.{\dynk \whtu0100{#1}\whtu1300{#2}\blku3000{#3}}
\def\eddgid#1.#2.#3.{\dynk \whtd0100{#1}\whtd1300{#2}\blkd3000{#3}}
\def\eddgiiu#1.#2.#3.{\dynk  \whtu0300{#1}\blku3100{#2}\blku1000{#3}}
\def\eddgiid#1.#2.#3.{\dynk  \whtd0300{#1}\blkd3100{#2}\blkd1000{#3}}

\def\ddfu#1.#2.#3.#4.{\dynk \whtu0100{#1}\whtu1200{#2}\blku2100{#3}\blku1000{#4}
}
\def\ddfd#1.#2.#3.#4.{\dynk \whtd0100{#1}\whtd1200{#2}\blkd2100{#3}\blkd1000{#4}
}
 
\def\eddfiu#1.#2.#3.#4.#5.{\dynk \whtu0100{#1}\whtu1100{#2}\whtu1200{#3}\blku210
0{#4}\blku1000{#5}}
\def\eddfid#1.#2.#3.#4.#5.{\dynk \whtd0100{#1}\whtd1100{#2}\whtd1200{#3}\blkd210
0{#4}\blkd1000{#5}}
\def\eddfiiu#1.#2.#3.#4.#5.{\dynk \whtu0100{#1}\whtu1200{#2}\blku2100{#3}\blku11
00{#4}\blku1000{#5}}
\def\eddfiid#1.#2.#3.#4.#5.{\dynk \whtd0100{#1}\whtd1200{#2}\blkd2100{#3}\blkd11
00{#4}\blkd1000{#5}}

\def\ddanu#1.#2.#3.#4.#5.{\dynk \whtu0100{#1}\whtu1100{#2}\whtu1100{#3}\cdots
                           \whtu1100{#4}\whtu1000{#5}}
\def\ddand#1.#2.#3.#4.#5.{\dynk \whtd0100{#1}\whtd1100{#2}\whtd1100{#3}\cdots
                           \whtd1100{#4}\whtd1000{#5}}

  \def\ddanm#1.#2.#3.#4.#5.{\dynk \whtu0100{#1}\whtu1100{#2}\whtu1100{#3}
                           \whtu1100{#4}\whtu1000{#5}} 
                           
    \def\ddnigel#1.#2.#3.{\dynk \whtu0100{#1}\whtu1100{#2}\whtu1100{#3}}
                                                       
       \def\ddanmm#1.#2.#3.#4.{\dynk \whtu0100{#1}\whtu1100{#2}
                           \whtu1100{#3}\whtu1000{#4}} 
                           
       \def\ddanmmm#1.#2.#3.{\dynk \whtu0100{#1}\whtu1100{#2}
                           \whtu1000{#3}}

 \def\ddans#1.#2.#3.#4.#5.{\dynk \blkd0100{#1}\whtd1100{#2}\blkd1100{#3}
                           \blkd1100{#4}\whtd1000{#5}}
\def\ddandto#1.#2.#3.#4.#5.{\dynk \blkd0100{#1}\whtd1100{#2}\blkd1100{#3}\cdots
                           \whtd1100{#4}\blkd1000{#5}}
 
\def\eddanu#1.#2.#3.#4.#5.{\dynk \whtu0100{#1}\whtu1100{#2}%
                           \up1{\whtr0000{#3}}\cdots\whtu1100{#4}\whtu1000{#5}}
\def\eddand#1.#2.#3.#4.#5.{\dynk \whtd0100{#1}\whtd1100{#2}%
                           \up1{\whtr0000{#3}}\cdots\whtd1100{#4}\whtd1000{#5}}
\def\eddaid#1.#2.{\dynk\whtd0400{#1}\hskip30pt\whtd4000{#2}}

\def\eddanid#1.#2.#3.#4.#5.{\dynk \whtd0200{#1}\whtd2100{#2}%
                           \whtd1100{#3}\cdots\whtd1200{#4}\blkd2000{#5}}
\def\eddaniu#1.#2.#3.#4.#5.{\dynk \whtu0200{#1}\whtu2100{#2}%
                           \whtu1100{#3}\cdots\whtu1200{#4}\blku2000{#5}}

\def\eddaniid#1.#2.#3.#4.#5.#6.{\hbox{$\vcenter{\hbox
         {\dynk\hbox{$ \lbde{#1}{#2}{#3}\blkd1100{#4}\cdots%
          \blkd1200{#5}\whtd2000{#6} $}} }$}}
\def\eddaniiu#1.#2.#3.#4.#5.#6.{\hbox{$\vcenter{\hbox
         {\dynk\hbox{$ \lbde{#1}{#2}{#3}\blku1100{#4}\cdots%
          \blku1200{#5}\whtu2000{#6} $}} }$}}

\def\eddaiiid#1.#2.{\dynk\blkd0400{#1}\hskip30pt\whtd4000{#2}}

\def\ddbnu#1.#2.#3.#4.#5.{\dynk \whtu0100{#1}\whtu1100{#2}\whtu1100{#3}\cdots
                           \whtu1200{#4}\blku2000{#5}}
\def\ddbnd#1.#2.#3.#4.#5.{\dynk \whtd0100{#1}\whtd1100{#2}\whtd1100{#3}\cdots
                           \whtd1200{#4}\blkd2000{#5}}
 
\def\eddbnu#1.#2.#3.#4.#5.#6.{\dynk \lde{#1}{#2}{#3}\whtu1100{#4}\cdots
                           \whtu1200{#5}\blku2000{#6}}
\def\eddbnd#1.#2.#3.#4.#5.#6.{\dynk \lde{#1}{#2}{#3}\whtd1100{#4}\cdots
                           \whtd1200{#5}\blkd2000{#6}}

\def\ddcnu#1.#2.#3.#4.#5.{\dynk \blku0100{#1}\blku1100{#2}\blku1100{#3}\cdots
                           \blku1200{#4}\whtu2000{#5}}
\def\ddcnd#1.#2.#3.#4.#5.{\dynk \blkd0100{#1}\blkd1100{#2}\blkd1100{#3}\cdots
                           \blkd1200{#4}\whtd2000{#5}}
 
\def\eddcnu#1.#2.#3.#4.#5.#6.{\dynk \whtu0200{#1}\blku2100{#2}\blku1100{#3}
       \blku1100{#4}\cdots
                           \blku1200{#5}\whtu2000{#6}}
\def\eddcnd#1.#2.#3.#4.#5.{\dynk \whtd0200{#1}\blkd2100{#2}\blkd1100{#3}
       \cdots \blkd1200{#4}\whtd2000{#5}}

\def\dddnu#1.#2.#3.#4.#5.#6{\hbox{$\vcenter{\hbox
         {\dynk\hbox{$ \whtu0100{#1}\whtu1100{#2}\cdots%
          \whtu1100{#3}\rde{#4}{#5}{#6} $}}  }$}}
\def\dddnd#1.#2.#3.#4.#5.#6.{\hbox{$\vcenter{\hbox
         {\dynk\hbox{$ \whtd0100{#1}\whtd1100{#2}\cdots%
          \whtd1100{#3}\rde{#4}{#5}{#6} $}} }$}}
\def\dddndte#1.#2.#3.#4.#5.#6.{\hbox{$\vcenter{\hbox
         {\dynk\hbox{$ \blkd0100{#1}\whtd1100{#2}\cdots%
          \blkd1100{#3}\rdet{#4}{#5}{#6} $}} }$}}
          
\def\dddndto#1.#2.#3.#4.#5.#6.{\hbox{$\vcenter{\hbox
         {\dynk\hbox{$ \whtd0100{#1}\blkd1100{#2}\cdots%
          \blkd1100{#3}\rdet{#4}{#5}{#6} $}} }$}}
\def\dddiv#1.#2.#3.#4.{\hbox{$\vcenter{\hbox
         {\dynk\hbox{$ \whtu0100{#1}\rde{#2}{#3}{#4}
              $}}  }$}}

\def\edddnu#1.#2.#3.#4.#5.#6.#7.#8.{\hbox{$\vcenter{\hbox
         {\dynk\hbox{$ \lde{#1}{#2}{#3}\whtu1100{#4}\cdots%
          \whtu1100{#5}\rde{#6}{#7}{#8} $}}  }$}}
\def\edddnd#1.#2.#3.#4.#5.#6.#7.#8.{\hbox{$\vcenter{\hbox
         {\dynk\hbox{$ \lde{#1}{#2}{#3}\whtd1100{#4}\cdots%
          \whtd1100{#5}\rde{#6}{#7}{#8} $}} }$}}

\def\edddniid#1.#2.#3.#4.#5.{\hbox{$\vcenter{\hbox
         {\dynk\hbox{$ \blkd0200{#1}\whtd2100{#2}\whtd1100{#3}\cdots%
          \whtd1200{#4}\blkd2000{#5} $}} }$}}
\def\edddniiu#1.#2.#3.#4.#5.{\hbox{$\vcenter{\hbox
         {\dynk\hbox{$ \blku0200{#1}\whtu2100{#2}\whtu1100{#3}\cdots%
          \whtu1200{#4}\blku2000{#5} $}} }$}}

\def\ddei#1.#2.#3.#4.#5.#6.{\hbox{$\vcenter{\hbox
       {\dynk \whtd0100{#1}\whtd1100{#3}%
       \up1{\whtr0001{#2}}\whtd1110{#4}\whtd1100{#5}\whtd1000{#6}} }$}}
\def\ddeit#1.#2.#3.#4.#5.#6.{\hbox{$\vcenter{\hbox
       {\dynk \whtd0100{#1}\blkd1100{#3}%
       \up1{\blkr0001{#2}}\whtd1110{#4}\blkd1100{#5}\whtd1000{#6}} }$}}
 
\def\eddei#1.#2.#3.#4.#5.#6.#7.{\hbox{$\vcenter{\hbox
       {\dynk \whtd0100{#1}\whtd1100{#3}%
       \up1{\whtr0011{#2}}\up2{\whtr0001{#7}}\whtd1110{#4}\whtd1100{#5}%
       \whtd1000{#6}} }$}}

\def\ddeii#1.#2.#3.#4.#5.#6.#7.{\hbox{$\vcenter{\hbox
       {\dynk \whtd0100{#1}\whtd1100{#3}%
       \up1{\whtr0001{#2}}\whtd1110{#4}\whtd1100{#5}\whtd1100{#6}%
       \whtd1000{#7}} }$}}
\def\ddeiit#1.#2.#3.#4.#5.#6.#7.{\hbox{$\vcenter{\hbox
       {\dynk \whtd0100{#1}\blkd1100{#3}%
       \up1{\blkr0001{#2}}\whtd1110{#4}\blkd1100{#5}\whtd1100{#6}%
       \blkd1000{#7}} }$}}
 
\def\eddeii#1.#2.#3.#4.#5.#6.#7.#8.{\hbox{$\vcenter{\hbox
       {\dynk \whtd0100{#8}\whtd1100{#1}\whtd1100{#3}%
       \up1{\whtr0001{#2}}\whtd1110{#4}\whtd1100{#5}\whtd1100{#6}%
       \whtd1000{#7}} }$}}

\def\ddeiii#1.#2.#3.#4.#5.#6.#7.#8.{\hbox{$\vcenter{\hbox
       {\dynk \whtd0100{#1}\whtd1100{#3}%
       \up1{\whtr0001{#2}}\whtd1110{#4}\whtd1100{#5}\whtd1100{#6}%
       \whtd1100{#7}\whtd1000{#8}} }$}}
\def\ddeiiit#1.#2.#3.#4.#5.#6.#7.#8.{\hbox{$\vcenter{\hbox
       {\dynk \whtd0100{#1}\blkd1100{#3}%
       \up1{\blkr0001{#2}}\whtd1110{#4}\blkd1100{#5}\whtd1100{#6}%
       \blkd1100{#7}\whtd1000{#8}} }$}}
 
\def\eddeiii#1.#2.#3.#4.#5.#6.#7.#8.#9.{\hbox{$\vcenter{\hbox
       {\dynk \whtd0100{#1}\whtd1100{#3}%
       \up1{\whtr0001{#2}}\whtd1110{#4}\whtd1100{#5}\whtd1100{#6}%
       \whtd1100{#7}\whtd1100{#8}\whtd1000{#9}} }$}}

\newcommand{\DIV}{\mathbf{D}_{\mathbf 4}}
\newcommand{\R}{\mathbf{R}}
\newcommand{\C}{\mathbf{C}}

\newcommand{\Z}{\mathbf{Z}}

\newcommand{\PP}{{\mathbf {\rm P}}}

\newcommand{\Ak}{\mathbf{A}_{\mathbf k}}
\newcommand{\Dk}{\mathbf{D}_{\mathbf k}}
\newcommand{\Ek}{\mathbf{E}_{\mathbf k}}
\newcommand{\EVI}{\mathbf{E}_{\mathbf 6}}
\newcommand{\EVII}{\mathbf{E}_{\mathbf 7}}
\newcommand{\EVIII}{\mathbf{E}_{\mathbf 8}}

\begin{document}
\title{The central sphere of an ALE space}
 \author{Nigel Hitchin\\[5pt]}
 \maketitle
 \centerline{\it {Dedicated to the memory of Sir Michael Atiyah}}
  \section{Introduction}
  Under the supervision of Michael Atiyah, Peter Kronheimer in his  thesis \cite{PB1},\cite{PB2} gave a beautiful construction of hyperk\"ahler 4-manifolds which are asymptotically locally Euclidean (ALE),  that is they have an infinity which is modelled on $\C^2\backslash \{0\}/\Gamma$ for a finite subgroup $\Gamma\subset SU(2)$. This yielded existence, regularity and moduli following naturally from properties of the regular representation of $\Gamma$.  The manifold is produced using the hyperk\"ahler quotient construction by the action of a Lie group $G$ on a quaternionic vector space.  To obtain the metric one considers the zero set of the moment map ($3\dim G$ quadratic equations) takes the induced Euclidean metric and then the quotient metric after acting by $G$. This is some way from being explicit and the aim of this paper is to produce something a little more concrete.
  
  In the case where $\Gamma$ is a cyclic group, there always has been an explicit form introduced by Gibbons and Hawking \cite{GH} in 1978, but the author's attack during the same period on the dihedral group  \cite{CH} yielded rather unmanageable expressions. In this article we specialize in two ways. Firstly we consider a choice of moment map which allows a circle action on the 4-manifold.  In this case, among the complex structures of the hyperk\"ahler family one is the resolution of the singularity $\C^2/\Gamma$. The configuration of rational curves which replaces the singularity has the structure of a Dynkin diagram and there is a distinguished vertex  (the trivalent one except for the type $\Ak$) which represents a curve which is pointwise fixed by the circle action. It is the metric on this central sphere which we aim to calculate. In fact, thanks to \cite{F},\cite{K}, this metric uniquely determines the ALE metric, and we may regard this as boundary data. 
  
  A hyperk\"ahler manifold has complex structures $I,J,K$ associated to the quaternions. In our case complex structure $I$ is the resolution of the singularity and the sphere is a holomorphic rational curve. In complex structure $J$, the manifold is an affine surface in $\C^3$ and the sphere is a compact component of a real form. We shall describe the metric in terms of two objects: a symplectic structure and a conformal structure. The symplectic structure is just the canonical 2-form on the real affine surface. The conformal structure will be described by the restriction  to the real form of a meromorphic function on the complex affine surface. 
  
  The method we use is Penrose's twistor space, or nonlinear graviton, construction  which was applied  in  \cite{H1} to the $\Ak$ case. The twistor space itself is defined by the versal deformation of the singularity together with  its simultaneous resolution using the original approach of Brieskorn, Tyurina and others \cite{B},\cite{Ty},\cite{Pink}. The resolution for $\Ak$ and $\Dk$ gives a concrete meromorphic function but for $\Ek$, $k=6,7,8$, namely when $\Gamma$ is the binary tetrahedral, octahedral or icosahedral group, we need a less direct description. The transition occurs by focusing on the first occasion where a trivalent vertex appears in the Dynkin diagram --  the case of $\DIV$. The affine surface here compactifies to  a nonsingular projective cubic and we relate the explicit meromorphic function to some of the 27 lines on it.

 By blowing up points on this cubic surface we find compactifications in the other cases and describe the meromorphic function in terms of  a particular configuration of exceptional curves. The function is the unique one, up to scaling and adding a constant, with a pole on this divisor. The reader may feel that this is somewhat less than concrete, but explicitness can be bought at too high a price. In particular the explicit equations for the affine surfaces involve the invariant polynomials for the exceptional groups which we have avoided. The reader may refer to the formulas in \cite{KM} to see why. 
 
 While the complex meromorphic function restricted to a real sphere defines the complex structure we may also regard the complex affine surface as possessing a complexification of this conformal structure, which consist of two families of null curves, related by complex conjugation.  Each one of these,  we shall see, is defined by a pencil of rational curves, the configuration being a reducible member. 
  
  In recent years more examples of hyperk\"ahler manifolds have attracted attention and at various points of the paper we draw comparisons with the moduli space of Higgs bundles on a Riemann surface where the analogue of the sphere studied here is the moduli space of stable bundles in complex structure $I$ or the flat unitary connections in complex structure $J$. 
  
  The author thanks Andy Hanson, whose persistent questions led to this study, and EPSRC for support. It goes without saying that Michael Atiyah introduced me, along with much more, to the world of algebraic geometry and its applications. 
  \section{Kronheimer's construction}
  We outline here  the hyperk\"ahler quotient construction and Torelli theorem \cite{PB1}, \cite{PB2} of ALE spaces which for us asserts the existence of the metric and description of moduli. 
  
  We take $\Gamma\subset SU(2)$, a finite subgroup, then $R=L^2(\Gamma)$ is a Hermitian vector space and the antilinear involution $A\mapsto A^*$ makes $\End(R)$  a real vector space. 
  Now $\C^2$ is a quaternionic vector space in the sense of having an antilinear map $J$ such that $J^2=-1$. The tensor product of a quaternionic and real vector space is quaternionic hence $\C^2\otimes \End(R)$  is quaternionic as is the fixed point set $V=(\C^2\otimes \End(R))^{\Gamma}$. This vector space is a flat hyperk\"ahler manifold with K\"ahler forms $\omega_1,\omega_2,\omega_3$.
  
  Let $U(R)$ be the group of unitary transformations of $R$ and $U(R)^{\Gamma}$ the subgroup commuting with the action of $\Gamma$. This is a product of unitary groups $U(n_i)$, $0\le i\le k$ one for each irreducible representation of $\Gamma$, of dimension $n_i$. It acts on $V$ preserving the flat hyperk\"ahler structure but because of the adjoint action on $\End(R)$ the scalars $T$ act trivially, so $G= U(R)^{\Gamma}/T$ acts on $V$. To take a hyperk\"ahler quotient we need to define a hyperk\"ahler moment map -- a triple of moment maps $\mu=(\mu_1,\mu_2,\mu_3)$ for the three symplectic forms $\omega_1,\omega_2,\omega_3$ -- and because $G$ has a centre (defined by the scalars in each $U(n_i)$) there is a choice, which will give parameters for the metric.  
  
  The McKay correspondence  associates to each  factor of $U(R)^{\Gamma}$ a vertex of an extended Dynkin diagram of type $\Ak,\Dk,\EVI,\EVII,\EVIII$. The vertices of the ordinary diagram correspond to simple roots $\theta_i$ and the  trivial representation is the extra vertex associated to  $\theta_0=-\sum_1^kn_i\theta_i$. Mapping a certain scalar $i\pi_i$ in $\lie{u}(n_i)$ to $n_i\theta_i$ identifies the Lie algebra $\lie {z}$ of the centre with the Cartan subalgebra $\lie{h}$ of the root system and the  moment map is then of the form $\mu: V\rightarrow \lie{h}\otimes \R^3$. The ALE space  is constructed as the  hyperk\"ahler quotient $M = \mu^{-1}(a)/G$ so long as $a$ does not lie in $H\otimes \R^3$ where  $H$ is a root hyperplane. 
  
  Given the ALE space $M$, the  parameter $a$ is recovered from the cohomology classes of the three K\"ahler forms  of the complex structures $I,J,K$  -- the intersection pairing on $H_2(M,\Z)$ is given by the Cartan matrix. 
  This becomes more evident when we restrict to the case where we have a circle action.  
  
  Here we take $a\in \lie{h}\subset \lie{h}\otimes \R^3$ (in a Weyl chamber to avoid the root hyperplanes). Equivalently  the moment maps $\mu_2,\mu_3$  vanish. The scalar action of $e^{i\theta}$ on the $\C^2$ factor in $V=(\C^2\otimes \End(R))^{\Gamma}$ preserves $\omega_1$ but acts as $e^{2i\theta}$ on $\omega_2+i\omega_3$.  It then preserves $\mu_2=0=\mu_3$ and descends to the hyperk\"ahler quotient as an isometry acting  on the induced K\"ahler forms in the same way. Note however that $-1$ is contained in $\Gamma$ except for $\Ak$ where $k$ is even, and so the effective circle action is the quotient by $\pm 1$. We shall regard this as the action, where the holomorphic 2-form transforms by $e^{i\theta}$. In this case the complex structure $I$ gives $M$ the structure of the  resolution of the singular space $\C^2/\Gamma$ and the origin is replaced by a configuration of rational curves whose intersections are given by the Dynkin diagram.  
  
  Apart from the $\Ak$ series, which we shall discuss separately, there is a single trivalent vertex in each of these. Since a circle action on the sphere cannot have just 3 fixed points that rational curve is pointwise fixed. It is a totally geodesic surface in $M$, which we call the {\it central sphere} $S_0$.  On the other spheres $S_i$ in the resolution the action is a rotation with two fixed points.

\vskip .5cm
\centerline{\dddnu. . .{\mathbf S}_0. . .\hskip .5cm\ddei . . .{\mathbf S}_0. . . }

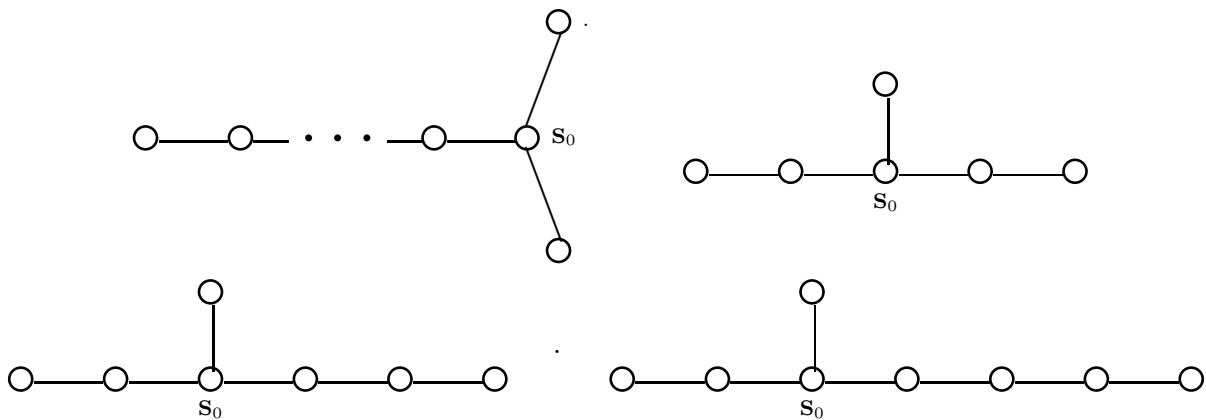
\centerline{\ddeii . . .{\mathbf S}_0. . . . . \ddeiii . . .{\mathbf S}_0. . . . . .}

\begin{rmks}

\noindent 1. The circle action preserves $\omega_1$ and has a moment map $f$ which is a proper Morse function with critical points at the fixed points. From the action on the tangent space we can see that the central sphere is a minimum. To obtain the other critical values note that the other spheres are surfaces of revolution and the restriction of $\omega_1$ is $d\theta\wedge df$, so integrating, or equivalently evaluating the cohomology class, gives $2\pi (f(q)-f(p))$ where $p$ and $q$ are the fixed points. If we choose $f=0$ on $S_0$ then this can be calculated in terms of the data $a\in \lie{h}$ defining the metric and the adjacency of the vertices. Note that since each sphere is holomorphic and $\omega_1$ is a K\"ahler form the integral is positive.

\noindent 2. In \cite{F},\cite{K} it is shown that any real analytic K\"ahler metric has a unique hyperk\"ahler extension to a circle-invariant metric on a neighbourhood of the zero section of the cotangent bundle. In our case, in complex structure $I$, $S_0$ is a rational curve of self-intersection $-2$ which means it can be collapsed to an ordinary double point and hence it is holomorphically equivalent to a neighbourhood of the zero section. This means that the ALE metric is uniquely determined by the metric on $S_0$. In particular, the constant curvature metric on the sphere determines the $\Ak$ ALE metric for $k=1$, the Eguchi-Hanson metric. The deformation-theoretic arguments in the above references can actually be replaced by a differential equation $g_{tt}=4\kappa g$ for a $t$-dependent family of metrics on the sphere where $\kappa$ is the Gaussian curvature \cite{H3}.

\noindent 3. In some respects the picture above is parallel to that in another well-studied hyperk\"ahler manifold, the moduli space ${\mathcal M}$ of Higgs bundles on a compact Riemann surface $\Sigma$ \cite{H2}. A Higgs bundle is  a holomorphic vector bundle $V$ together with a holomorphic section $\Phi$ of $\End V\otimes K$ satisfying a stability condition. This condition implies the existence of a Hermitian metric such that $F_A+[\Phi,\Phi^*]=0$ and the natural $L^2$ metric on the moduli space of solutions to this equation is hyperk\"ahler. Here we also have a circle action $\Phi\mapsto e^{i\theta}\Phi$ and a proper Morse function  $f=\Vert\Phi\Vert^2$. The absolute minimum of $f$ is $\Phi=0$ which is the moduli space of (semi)-stable vector bundles  and is the analogue of our sphere $S_0$. The Higgs bundles for which $\Phi$ is nilpotent forms the nilpotent cone, a topologically connected configuration of half-dimensional compact subvarieties  preserved by the action, analogous to the resolution of the singularity. A key difference is that the intersection form (when ${\mathcal M}$ is smooth) is zero \cite{Hein} whereas it is nondegenerate in the ALE case. 
\end{rmks}
\section{Multi-instanton metrics}\label{multi}
The multi-instanton metrics of Gibbons and Hawking \cite{GH} were the first examples of ALE spaces. The authors used a metric of the form 
$$g=V(dx_1^2+dx_2^2+dx_3^2)+V^{-1}(d\theta +A)^2$$
where $A$ defines a $U(1)$-connection on $\R^3$ and $V$ is a function on $\R^3$.  If $F_A=dA$ is the curvature and $dV=\ast F_A$ then $\omega_1=Vdx_2\wedge dx_3+dx_1\wedge (d\theta +A)$ and similar forms are closed and constitute the three K\"ahler forms of a hyperk\"ahler family. Since $dF_A=0$, the Ansatz reduces to the consideration of a single harmonic function $V$ on $\R^3$. 

The standard Euclidean metric on $\C^2$ can be put in this form by setting 
$$x_1=\frac{1}{2}(\vert z_1\vert^2-\vert z_2\vert^2)\qquad x_2+ix_3=z_1z_2$$
and then 
$$V=\frac{1}{\vert z_1\vert^2+\vert z_2\vert^2}=\frac{1}{2\sqrt{x_1^2+x_2^2+x_3^2}}.$$
In this case the principal $U(1)$-bundle (the Dirac monopole) over a  sphere in $\R^3$ is $S^3\subset \C^2$ or the Hopf bundle. Note that the circle action $(z_1,z_2)\mapsto (e^{i\theta}z_1, e^{i\theta}z_2)$ is a rotation by $2\theta$ about the $x_1$-axis in $\R^3$.

The multi-instanton solution consists of taking $k+1$ points ${\mathbf a_i}\in \R^3$ and setting 
$$V=\sum_{i=1}^{k+1}\frac{1}{2\vert {\mathbf x}-{\mathbf a_i}\vert}.$$
The $U(1)$-bundle now has degree $k+1$ over the sphere at infinity and the principal bundle is $S^3/\Z_{k+1}$. Infinity behaves then like $\C^2/\Gamma$ for the cyclic group. The apparent singularities at ${\mathbf a_i}\in \R^3$ are in fact smooth by comparison with the flat case. The global picture is of a 4-manifold $M$ with a $U(1)$-action preserving the hyperk\"ahler forms and having fixed points over  ${\mathbf a_1},\dots,{\mathbf a_{k+1}}$. The functions $(x_1,x_2,x_3)$ form the hyperk\"ahler moment map for this action.

\begin{rmk} The origin of this Ansatz lies in a solution to the Einstein-Maxwell equations \cite{Maj} in Lorentzian signature. Adapted to  Euclidean signature the energy-momentum tensor for a Maxwell field is the tensor product of the self-dual and anti-self-dual part so a self-dual 2-form gives a solution to the vacuum Einstein equations. 
\end{rmk}

For the circle-invariant solutions we are considering we take the points ${\mathbf a_i}$ to lie on the  $x_1$-axis and then rotation about that axis induces an isometric circle action generating a vector field $X$. This involves lifting the rotation by $2\theta$ on  $\R^3$ to the $U(1)$-bundle with connection form $d\theta+A$, commuting with the action. Such a lifting defines a vector field of the form 
$X=X_H+hY$, where $Y=\partial/\partial\theta$ is the vertical vector field on $M$ and 
 $X_H$ is the horizontal lift of $$2x_2\frac{\partial}{\partial  x_3}-2x_3\frac{\partial}{\partial x_2}$$ 
Since $X$ preserves the connection form $d\theta+A$  this gives
$$dh=2(x_2V_2+x_3V_3)dx_1-2x_2V_1dx_2-2x_3V_1dx_3.$$
It follows that, with ${\mathbf a_i}=(a_i,0,0)$,  
\begin{equation}
h=\sum_{i=1}^{k+1}\frac{x_1-a_i}{\vert {\mathbf x}-\mathbf {a_i}\vert}+c.
\label{eff}
\end{equation}
The constant $c$ reflects the fact that any lift of the action can be modified by incorporating  an action of $U(1)$. 

For any choice, the  zeros of $Y$ project to the points $a_i$ on the $x_1$-axis and these are zeros of $X$. Order the points $a_1<a_2<\cdots<a_{k+1}$  and the inverse image of the intervals $[a_i,a_{i+1}]$ form a chain of  spheres, holomorphic in complex structure $I$, which is the resolution of $\C^2/\Gamma$.

 A sphere which is pointwise fixed  by the action means choosing an interval for which $h=0$. When $k=2\ell-1$ is odd the  natural choice is to take $c=0$ and $x_1\in [a_{\ell},a_{\ell+1}]$, the middle interval. We shall take this to be our central sphere $S_0$.  
$${\ddanm..{\mathbf S}_0...}$$
Note that when $k$ is odd  the rotation action on $\R^3$ lifts to the ALE space $M$. 

The moment map $f$ for $X$ satisfies $i_X\omega_1=df$ and using   $\omega_1=Vdx_2\wedge dx_3+dx_1\wedge(d\theta+\alpha)$ this gives
$$f=\sum_{i=1}^{2\ell}\vert {\mathbf x}-{\mathbf a_i}\vert $$
as the Morse function whose minimum locus is $S_0$.

Over the real axis the $U(1)$ connection is trivial and so the metric on $S_0$ is the surface of revolution 
$$g=Vdx^2+V^{-1}d\theta^2$$
where $$V=\sum_1^{2\ell}\frac{1}{2\vert { x}-{ a_i}\vert}$$
and $x\in[a_{\ell},a_{\ell+1}]$.
Having an explicit form for the metric means that we can observe the resolution of the $\Ak$ singularity without using any algebraic geometry. This is exceptional however, and so we proceed in a different manner.
\section{The twistor construction}
 Penrose's nonlinear graviton construction converts the problem of finding a hyperk\"ahler manifold  into one of holomorphic geometry. In the 4-dimensional case considered here it requires first a complex 3-manifold $Z$, the {\it twistor space}, with a holomorphic fibration $\pi:Z\rightarrow \PP^1$. Additional data is a real structure: an antiholomorphic involution $\sigma$ compatible with the fibration and inducing the antipodal map $u\mapsto -1/\bar u$ on $\PP^1$ and a non-vanishing real section of the line bundle $K_Z(4)$, or better $\Lambda^2T^*_F(2)$, where  $T_F$ is the tangent bundle along the fibres tensored with $\pi^*{\mathcal O}(2)$.  This defines a holomorphic symplectic form on each fibre. 

Second, and more importantly, we require a family of holomorphic sections, {\it twistor lines}, with normal bundle ${\mathcal O}(1)\oplus {\mathcal O}(1)$. These belong to a complete complex 4-dimensional family with a conformal structure: two points are null-separated if the corresponding twistor lines intersect. The relative symplectic structure  
fixes a metric in the conformal class and a component of the space of real sections is the hyperk\"ahler manifold $M$.  Identifying the base $\PP^1$ with the 2-sphere of complex structures $S^2=\{aI+bJ+cK: a^2+b^2+c^2=1\}$ in the hyperk\"ahler family, the fibre of $\pi$ over a point is the complex manifold $(M,aI+bJ+cK)$. A circle action of the type we are considering gives a holomorphic action on $Z$, compatible with the fibration structure and acting on an affine coordinate $u$ on $\PP^1$ by $u\mapsto e^{i\theta}u$. 

Removing a line $L$ from $\PP^3$ and using $K_{\PP^3}\cong {\mathcal O}(-4)$ gives the flat $\C^2$ twistor space with  projective lines $\PP^1\subset \PP^3\backslash L$ being the twistor lines. Constructing the twistor space for an ALE metric uses the {\it simultaneous resolution} of the Kleinian singularity $\C^2/\Gamma$ \cite{Ty},\cite{Pink}. 

For $\Gamma\subset SU(2)$ the ring of $\Gamma$-invariant polynomials is generated by three functions $x,y,z$ which satisfy an algebraic relation which represents $\C^2/\Gamma$ as a surface in $\C^3$. For the cyclic group action $(z_1,z_2)\mapsto (\omega z_1,\omega^{-1}z_2)$  with $\omega^{k+1}=1$ we have $x=z_1^{k+1}, y=z_2^{k+1}, z=z_1z_2$ satisfying $xy=z^{k+1}$. The versal deformation of this singularity involves adding lower order terms to give a family of surfaces, in this example $xy=z^{k+1}+c_1z^k+\cdots+c_{k+1}$. The generic one is smooth but singular ones occur where the polynomial in $z$ has multiple roots. The simultaneous resolution of the whole family requires one  to parametrize the deformation by  the roots $a_i$ of this polynomial  and consider 
$$xy=\prod_{i=1}^{k+1}(z-a_i)=z^{k+1}+c_1z^k+\cdots+c_{k+1}$$
With $c_1=0$, the parameter space $\C^k=\{(a_1,\dots, a_{k+1}): \sum _i a_i=0\}$ is recognizable as the complex Cartan subalgebra ${\lie h}^c$ of type $\Ak$ and the coefficients $c_i$ in the equation of the surface as generating polynomials for the ring of invariants under the Weyl group. This is the general situation. 

The versal deformation defines  a subvariety $\tilde Y\subset \C^3\times {\lie h}^c$ and a simultaneous resolution is a smooth variety $Y\rightarrow {\lie h}^c$ with a projection $Y\rightarrow \tilde Y$ over ${\lie h}^c$ which fibre-by-fibre is a minimal resolution of the singular surfaces in the deformation. The original construction for these singularities  is a case-by-case treatment due to Brieskorn and Tyurina but there are other more uniform constructions such as the Slodowy slice of a subregular nilpotent.  It is more convenient here to use the original approach. 

The resolution has some fundamental properties:
\begin{itemize}
\item
The invariant polynomials $x,y,z$ are homogeneous of  degrees $p,q,r$ in $(z_1,z_2)$ and the scalar action of $\C^*$ on $\C^2$ induces an action on $Y$ commuting with the scalar action on ${\lie h}^c$.
\item
The holomorphic symplectic form $dz_1\wedge dz_2$ on $\C^2$ is $\Gamma$-invariant and defines one on the resolution of $\C^2/\Gamma$. This belongs to a relative  symplectic form on the fibres of $Y\rightarrow {\lie h}^c$ which on a smooth affine surface $f(x,y,z)=0$ over a point in ${\lie h}^c$ is a multiple of the standard form $dx\wedge dy/f_z$.
\item
When the degrees of $x,y,z$ are even (the only exception is $\Ak$ for $k$ even), the quaternionic structure $(z_1,z_2)\mapsto (\bar z_1,-\bar z_1)$ induces a real structure on $Y$ compatible with $\lie{h}^c=\lie{h}\otimes \C$.
\end{itemize}
These properties enable one to define the twistor space for Kronheimer's construction. The value of the moment map in $\lie{h}\otimes \R^3$ we take to be a real section of ${\lie h}^c(2)$ on $\PP^1$ -- quadratic polynomials in the affine parameter $u$ with values in ${\lie h}^c$. For $u\ne \infty$ this is a  map $h:\C\rightarrow \lie{h}^c$ and we define $Z_+$ to be the pullback $h^{*}Y$, a nonsingular 3-manifold resolving the singular fibres in the 1-dimensional family of deformations. 
The pull-back of the versal deformation defines a singular subvariety $\tilde  Z\subset {\mathcal O}(p)\oplus {\mathcal   O}(q)\oplus{\mathcal O}(r)\rightarrow \PP^1$ 
with the real   structure acting antipodally on $\PP^1$. The space $Z_+$ resolves the singularities for $u\ne \infty$ and applying the real structure we obtain $Z_-$ doing the same for $u\ne 0$. Patched together, they form the twistor space $Z\rightarrow \PP^1$  resolving $\tilde Z$.

 It will be sufficient for most purposes to work with just $Z_+$. So this is a fibration over $\C$, the fibre $Z_0$ is the resolution of $\C^2/\Gamma$ and a general fibre $Z_u$ over $ u\ne 0 \in \C$ is an affine surface in $\C^3$.
\begin{rmks} 

\noindent 1. Note that although $Z_+$ is quasi-projective, the twistor space   $Z$ is not, as the sphere $S_0$ in the fibre over $u=0$ has the opposite complex structure and hence orientation to the sphere over $u=\infty$, yet they are in the same cohomology class.  

\noindent 2. In the case of Higgs bundles, the corresponding space $Z_+$  can be identified (with $u$ replacing $\lambda$) with the moduli space of $\lambda$-connections -- holomorphic differential operators on sections of a vector bundle such that $\nabla(fs) =f\nabla s+\lambda s\otimes df$. If $\lambda \ne 0$ then $\lambda^{-1}\nabla$ is just a usual holomorphic flat connection. 
\end{rmks} 

\begin{ex} The twistor space for  $\Ak$ is defined by  quadratic polynomials $p_i(u)=b_iu^2+a_iu-\bar b_i$ and setting $\tilde Z$ to be  
$$xy=\prod_{i=1}^{k+1}(z-p_i(u)).$$
Finding the sections, the twistor lines,  means finding polynomial solutions  $x(u),y(u)$ and $z(u)$ of degrees $k+1,k+1,2 $ to this equation. This was carried out for generic sections in \cite{H1} by simple factorization where it gave an alternative expression for the metric in Section \ref{multi}, better adapted to the complex structure.
\end{ex}

We described above the twistor space for a general ALE space. This paper concerns the case where we have a circle action. Real sections of ${\mathcal O}(2)$ on $\PP^1$ which are invariant by the action are polynomials of the form $cu$ where $c$ is real. In this case we take $h(u)=au$ for $a\in {\lie h}$. Then the circle action on $Z_+$ is just the natural action on the simultaneous resolution. To describe $S_0$ we  only need the twistor lines which correspond to fixed points of the circle action: this means such a line is the closure of an orbit of the $\C^*$-action. 

If $a\in {\lie h}$ does not lie on any root plane then the fibre $Z_1$ over $u=1$ (or complex structure $J$ in the standard parametrization) is a nonsingular affine surface. The  circle  acts as $u\mapsto e^{i\theta}u$ and rotates the $J,K$-plane of complex structures in the hyperk\"ahler family.  So $\theta=\pi$ takes $J$ to $-J$. It follows that the fixed point set of the circle is also fixed by an antiholomorphic involution -- a real structure. Since $S_0$ is 2-dimensional it  must be a connected component of a real form of the affine surface. Hence a twistor line defining a point in $S_0$ is defined by taking a real point of the surface $Z_1$ and its orbit under the $\C^*$-action.
\begin{ex}
In the ${\mathbf A}_{\mathbf {2\ell-1}}$ case we have $p_i(u)=a_iu$ where $a_i$ is real and $a_i\ne a_j$ to define the twistor space. The real structure  is $(x,y,z)\mapsto (\bar y,\bar x, \bar z)$, so complex $x$ and real $z$ satisfying 
$$x\bar x=\prod_{i=1}^{2\ell}(z-a_i)$$
is a real point and the twistor line is $u\mapsto(xu^{\ell},\bar xu^{\ell}, zu)$. If the $a_i$ are ordered, then setting $u=1$ we see two noncompact components for $z\le a_1$ and $z\ge a_{2\ell}$ and disjoint spheres for $z\in (a_2,a_3), (a_4,a_5)$ etc. In complex structure $I$, the fibre $Z_0$, these are closures of $\C^*$-orbits and are fixed points of a {\it holomorphic} involution.
\end{ex} 
\begin{rmk} 
In the analogous case of Higgs bundles, complex structure $J$ is the moduli space of flat connections for a complex group $G^c$ and the real points correspond to the holonomy being in a real form of $G^c$. There is generally a single compact component corresponding to the compact real form. As in the case here, each component  in complex structure $I$ is holomorphic and $\C^*$-invariant which provides the opportunity to use Morse theory for the other real forms as in \cite{H2} and the many results of O.Garcia-Prada  and collaborators e.g. \cite{GP}. The most studied  problem there is to determine the connected components for a non-compact real form of $G^c$ by looking for a local minimum of the Morse function. 

In our case note that the action on the tangent space at a point of intersection of $S_0$ with another sphere $S_1$ is $(1,e^{i\theta})$ and at the other fixed point  on $S_1$ the action on its tangent space is 
$e^{-i\theta}$. Since the $I$-holomorphic symplectic  form is acted on by $e^{i\theta}$ the action on $S_2$ at this point must be $e^{2i\theta}$ and so on. In particular $\theta=i\pi$ gives complex conjugation in complex structure $J$ so the  spheres $S_2$, $S_4$ etc. along any branch of the Dynkin diagram are compact real components. Similarly, if the final sphere on a branch is not real, there is a noncompact real component intersecting it, where the point of contact is the minimum of the Morse function. Thus, for example, the $\DIV$ surface has a unique compact component and 3 noncompact ones, $\EVI$ has 3 compact and one noncompact, $\EVII$ 3 compact and 2 noncompact and $\EVIII$ 4 compact and one noncompact.

\end{rmk}
The real surface together with the real symplectic 2-form $\omega$ define $S_0$ as a symplectic manifold. To obtain the metric we now need the conformal structure. The twistor construction for a complex spacetime tells us that two points are null separated if the twistor lines intersect. In this case it means they have same limit as $u\rightarrow 0$ or $\infty$ which of course is in the resolution of the singularity. 

We observe that $Z_+$ is a  nonsingular 3-dimensional variety  with a $\C^*$-action and as such it has a Bialynicki-Birula decomposition into locally closed subvarieties, each of which is an affine bundle over a component of the fixed point set. The curve $S_0\subset Z_0$ is the unique 1-dimensional fixed point set and the associated subvariety is a dense open subset of $Z_+$, the total space of a  bundle of 2-dimensional affine spaces. The affine bundle over a point in $S_0$ is then a union of $\C^*$-orbits comprising a surface, which intersects $Z_1$ in a curve.
This, by the definition of conformal structure,  is a null curve in the affine surface $Z_1$. Varying the point in $S_0$ we have the foliation by null curves which, together with their complex conjugates, defines a  complex conformal structure invariant under conjugation.

 On the real surface, a conformal structure is the same as a complex structure and this is realized by projection onto the quotient space of the foliation in its complexification. We deduce then:

\begin{prp}  Let $M$ be  an ALE space defined by $a\in {\lie h}$.  Then,
\begin{itemize}
\item
the central sphere $S_0$ can be identified with a compact component of a real form of the affine surface defined by $a$ in the versal deformation. 
\item
the volume form is a multiple of the canonical 2-form on the affine surface, 
\item
the complex structure is defined by the projection onto the fixed rational curve provided by the Bialynicki-Birula decomposition of the 3-fold $Z_+$.
\end{itemize}
\end{prp} 
\section{The central spheres of type $\Ak$ and $\Dk$}
\subsection{The ${\mathbf A}_{\mathbf{2\ell-1}}$ case}\label{Acentral} 
As discussed in Section \ref{multi} we restrict to the odd case $k=2\ell-1$, but now pursue the twistor approach. This means using the explicit simultaneous resolution as produced in \cite{B},\cite{Ty} and in more detail in \cite{KM}. We start with the versal deformation 
$$\tilde Y=\{(x,y,z,a)\in \C^3\times \C^k: \sum_{i=1}^{2\ell}a_i=0\}$$
and define $\mu:\tilde Y\rightarrow (\PP^1)^{2\ell-1}$ by 
$\mu_j=[x,\prod_1^j(z-a_i)]$
in homogeneous coordinates. Then the simultaneous resolution $Y$ is the closure of the graph of $\mu$. 

Using affine coordinates $v_i$ on the $\PP^1$ factors,  $u \in\C^*$ acts as 
$$(u^{\ell}x,u^{\ell}y, u z, u^{\ell-1}v_1,u^{\ell-1}v_2,\dots,u^2 v_{\ell-1}, v_{\ell}, u^{-2}v_{\ell+1},\dots, u^{-\ell+1}v_{2\ell-1})$$
so the generic limit as $u\rightarrow 0$ is $(0,0,\dots, v_{\ell},\infty,\dots,\infty)$. This projective line therefore lies in the closure of the graph over the origin in $\C^3$, is $\C^*$-invariant  and is a component of the resolution.  Note that the ordering in $(\PP^1)^k$ is the ordering of the rational curves in the $\Ak$ Dynkin diagram. Circle actions modified by the other $U(1)$-action  give any factor as a limit.

We can describe the situation more analytically by saying that  the meromorphic function
$$v_{\ell}=\frac{x}{\prod_{i=1}^{\ell}(z-a_i)}$$
restricted to the real component of the surface is a holomorphic coordinate, or more globally an identification of $S_0$ with $\PP^1$. Using this we are in a position to calculate the metric on $S_0$. 

The real points are defined by $x=\bar y=re^{i\theta}$ and taking $z$ real where $r^2=\prod_1^{2\ell}(z-a_i)$ so taking the logarithmic derivative of $v_{\ell}$ a $(1,0)$-form is 
$$\frac{dr}{r}+id\theta-\sum_{i=1}^{\ell}\frac{dz}{(z-a_i)}.$$
But from $r^2=\prod_1^{2\ell}(z-a_i)$
$$2\frac{dr}{r}=\sum_{i=1}^{2\ell}\frac{dz}{(z-a_i)}$$
so the $(1,0)$-form is 
$$id\theta-\frac{1}{2}\sum_{i=1}^{\ell}\frac{dz}{(z-a_i)}+\frac{1}{2}\sum_{i=\ell+1}^{2\ell}\frac{dz}{(z-a_i)}.$$
If the $a_i$ are ordered as before this is non-zero for $z$ in an interval $(a_i,a_{i+1})$ only if it is the central one $i=\ell$ and then this expression is 
$$id\theta +\sum_1^{2\ell}\frac{dz}{2\vert { x}-{ a_i}\vert}=id\theta+Vdz$$ as in the Gibbons-Hawking form.

The symplectic form for the metric $Vdz^2+V^{-1}d\theta^2$ is $dz\wedge d\theta$ and the canonical form on the surface $xy=\prod_1^{2\ell}(z-a_i)$ is $dx\wedge dz/x$. Since $x$ here is not a real coordinate this is not a real 2-form, in fact  on the real component it is $dx\wedge dz/ x=-idz\wedge d\theta$. So $idx\wedge dz/x$  recovers the standard form for the multi-instanton metric. 
\subsection{The  $\Dk$ case}
 This is the case where $\Gamma$ is the binary dihedral group of order $4k$. The extra symmetry in the dihedral group is $(z_1,z_2)\mapsto (z_2,-z_1)$ so $z=z_1^2z_2^2$ is invariant. The other generators are  $x=(z_1^{2k}+z_2^{2k})/2$ and 
 $y=(z_1^{2k}-z_2^{2k}) z_1z_2/2$ and the relation is given by 
 $x^2-zy^2=-z^{k+1}$. The versal deformation here is 
$$x^2-zy^2=-\frac{1}{z}\left(\prod_{i=1}^k(z+a_i^2)-\prod_{i=1}^ka_i^2\right)+2y\prod_{i=1}^k a_i$$
 with $\C^*$-action $(x,y,z,a)\mapsto (u^kx,u^{k-1}y, u^2z, ua)$ 
 and the real structure $(x,y,z)\mapsto (\bar x,\bar y,\bar z)$.  We need $\pm a_i\pm a_j\ne 0$ to keep away from the root planes of $\Dk$. 
 
 To consider the real forms, rewrite the equation  as  \begin{equation}
 -zx^2+(zy+\prod_i a_i)^2=\prod_i(z+a_i^2) 
  \label{alt}
 \end{equation}
 and note that for real values, fixed $z>0$ gives a hyperbola and $z<0$  gives a circle so long as the right hand side is positive. Ordering $a_1^2>a_2^2\cdots >a_k^2$ we then have compact components of the real form for $-z\in (a_2^2,a_3^2), (a_4^2,a_5^2). ...$. 
 
 Setting $z=-s^2$ for $s\in (a_i,a_{i+1})$ and $p=\prod_ia_i$, define $R^2=\prod_i(a_i^2-s^2)$ then $x=R\cos\theta/s, y=(-p+R\sin\theta)/s^2$  gives $(R,\theta)$ as standard coordinates on a sphere.
 Moreover, the symplectic form $\omega=dx\wedge dz/f_y=dx\wedge dz/(-2zy-2p)=ds\wedge d\theta$ so evaluating on the sphere we obtain $\pm 2\pi(a_{i+1}-a_i)$, the sign depending on orientation. 
 
The standard labelling of the vertices of the Dynkin diagram by simple roots is the following:  \vskip .5cm
\centerline{\dddnu{x_k-x_{k-1}}. .{x_4-x_3}.{x_3-x_2}.{-x_1+x_2}.{x_1+x_2}.}
\vskip .5cm
If $a\in \lie{h}$ lies in this positive Weyl chamber then $a_k>a_{k-1}>\cdots>a_2>a_1$ and $a_1+a_2>0$. So $a_2^2-a_1^2=(a_2-a_1)(a_2+a_1)>0$ and $a_i>0$ for $i\ge 2$ so the $a_i^2$ are ordered by magnitude as above. Our evaluation of $\omega$ on the real components then shows that  $-z\in (a_2^2,a_3^2)$ on $S_0$.

To produce the simultaneous resolution we  follow Tyurina \cite{Ty}. 
 First, considering the even and odd terms, we have polynomials $P,Q$  satisfying 
$$\prod_{i=1}^k(v+a_i)=v^k+\sigma_1v^{k-1}+\sigma_2v^{k-2}+\dots + p= vP_1(v^2)+Q_1(v^2).$$
Note that $Q(0)=a_1a_2\dots a_k=p$. 
Then
$-v^2P_1^2(v^2)+Q_1^2(v^2)=\prod_{i}(-v^2+a_i^2)$
or, setting  $P(z)=P_1(-z), Q(z)=Q_1(-z)$,
$$zP(z)^2+Q^2(z)=\prod_{i=1}^k(z+a_i^2).$$
Since $Q(0)=p$ we define $S$ by $Q(z)-p=zS(z)$. 
Then the equation of the surface is 
$$x^2-zy^2=-\frac{1}{z}\left(zP^2+Q^2-p^2\right)+2yp=-P^2-2pS-zS^2+2py$$
or
$(x+iP)(x-iP) =(y-S)(z(y+S)+2p).$

Define $t=(x+iP)/(y-S)$ then the  equation becomes 
$t(x-iP)=(z(y+S)+2p)$ or, eliminating $x$, 
\begin{equation}
(t^2-z)(y-S)=2zS+2p+2iPt=2(Q+iPt).
\label{blow}
\end{equation}
From the definition of $P$ and $Q$, $Q_1(-t^2)+itP_1(-t^2) =\prod_1^k(it+a_i)$
so that 
$$Q(z)+itP(z)-\prod_{i=1}^k(it+a_i)$$
is divisible by $z-t^2$, with factor $G(z)$. Then equation (\ref{blow})  becomes 
$$(z-t^2)(y+2G-S)=-2\prod_{i=1}^k(a_i+it)$$
which has the same format as ${\mathbf A}_{\mathbf{k-1}}$.

As in the $\Ak$ case, the resolution is defined by the closure of a subvariety in $\C^3\times (\PP^1)^k$ where $t$ takes values in the first factor (the $x_1+x_2$ rational curve) and the other ones are $[z-t^2,\prod_1^j(a_i+it)]$. The action of $u\in \C^*$ is $t\mapsto ut$, and so is  trivial if $j=2$. It follows that  the conformal structure is defined by the meromorphic function 
\begin{equation}
\frac{(a_1+it)(a_2+it)}{z-t^2}.
\label{mero}
\end{equation}
In principle, with the  coordinates $(R,\theta)$ above, this gives an explicit form of the metric on $S_0$ but we have lost the geometry. We next look at the meromorphic function (\ref{mero})  for the smallest value $k=4$ from a different viewpoint, for this will provide a link to the less accessible cases of $\Ek$. 
\section{The $\DIV$ cubic surface} 
When $k=4$ the affine surface is cubic. Its projective completion in homogeneous coordinates  is 
$$wx^2-zy^2=-\frac{1}{z}\left(\prod_{i=1}^4(z+wa_i^2)-w^4\prod_{i=1}^4a_i^2\right)+2w^2yp$$
and can easily be checked to be nonsingular. The plane at infinity $w=0$ intersects the surface in $zy^2=z^3$, three real lines $z=0, y=z, y=-z$. Call these $E,F,G$. They intersect in a  common point $(1,0,0,0)$.

A nonsingular projective cubic surface famously has 27 lines, moreover a set of 6 disjoint lines $E_i$ can be blown down to give 6 points $e_1,\dots, e_6\in \PP^2$. The further 21 are then provided by the proper transform $E_{ij}$ of the 15 lines in the plane joining the  points $e_i$ to $e_j$ and the 6 conics $C_i$ through five of the points, missing $e_i$. 

Consider the form of the equation (\ref{alt}) 
$$ -zx^2+(zy+p)^2=\prod_{i=1}^4(z+a_i^2).$$ 
If $z=-a_i^2$ then $a_i^2x^2+(p-ya_i^2)^2=(ia_ix+p-ya_i^2)(-ia_ix+p-ya_i^2)=0$, so the intersection of the planes $z+a_i^2=0$ and $ia_ix+p-ya_i^2=0$ is a line $E_j$ lying on the cubic. Its conjugate $\bar E_i$ is given by the other factor $-ia_ix+p-ya_i^2=0$. The plane $z+a_i^2=0$ intersects the cubic in these two lines together with the line $E$, given by $w=z=0$. 

We now have five disjoint lines $E_1,E_2,E_3,E_4, F$ which can be blown down to points $e_1,e_2,e_3,e_4, f$ in $\PP^2$ with a point blown up. None of the points meets the exceptional curve here for otherwise we would have a $-2$ curve in a smooth cubic, so there is a further line, call it $X$, in the cubic disjoint from these five, and blowing down all six gives $\PP^2$. 

The intersection properties of $X$ show that it meets $\bar E_1,\bar E_2,\bar E_3,\bar E_4$ and $F$. Using decomposable 2-forms this is enough to determine its equation which is $$x+i\sigma_1y-i\sigma_3=0,\qquad y-z+\sigma_2 =0$$
where $\sigma_i$ as usual are the elementary symmetric functions of the $a_i$. Tyurina's parameter $t$ is given in this case by
\begin{equation}
t=\frac{x+i\sigma_1z-i\sigma_3}{y-z+\sigma_2}
\label{teq}
\end{equation}
so geometrically $t= a$ is a plane through the line $X$, which meets the line $E$ at infinity in the point $(a,1,0,0)$.

The resulting picture of the cubic realized as the blow up of $\PP^2$ is Figure \ref{pic}. 
\vskip .25cm
\begin{figure}
\hskip 4cm\begin{tikzpicture}
\draw[] (4,3) ellipse (3cm and 1.5cm);
\draw[] (1,0) -- (1,6);
\draw[] (1,0) -- (2.5,6);
\draw[] (1,0) -- (4,6);
\draw[] (1,0) -- (5.7,6);
\draw[] (1,0) -- (7.8,6);
\node at (1,6.5) {${G}$};
\node at (2.5,6.5) {$\bar{ E}_1$};
\node at (4,6.5) {$\bar{ E}_2$};
\node at (5.7,6.5) {$\bar{ E}_3$};
\node at (7.8,6.5) {$\bar{ E}_4$};
\node at (0.75,3) {${ f}$};
\node at (0.75,0) {${ x}$};
\node at (4,1.25) {${ E}$};
\node at (1.75,4.35) {${ e}_1$};
\node at (3,4.7) {${ e}_2$};
\node at (4.25,4.75) {${ e}_3$};
\node at (5.7,4.5) {${ e}_4$};
\end{tikzpicture}
\caption{The $\DIV$ cubic surface}
\label{pic}
\end{figure}

Note that $E$ appears as the conic which misses the point $x\in \PP^2$. Let $C_1,\dots, C_4$ be lines defined by the conics which miss $e_1,\dots, e_4$. We shall give a more geometric description of the function 
$$f=\frac{(a_1+it)(a_2+it)}{z-t^2}.$$

\begin{prp} \label{D4div} Let $V$ be a projective nonsingular cubic surface in the $\DIV$ versal family and let $\hat V$ denote the surface obtained by blowing up the common point of intersection of the lines $E,F,G$. Let $\hat D$ denote the divisor which is the proper transform of $E+\bar E_3+\bar E_4$. 

Then the linear system $\vert \hat D\vert$ is a pencil of rational curves defining the meromorphic function $f$. More concretely, up to scalar multiplication and addition of a constant, $f$ is the unique meromorphic function with polar divisor $\hat D$.
\end{prp}

\begin{prf}

We calculate on the cubic surface $V$.
The numerator in the meromorphic function  has a factor $a_1+it$ and for $k=4$, $a_1+it=0$ is a plane which intersects the line at infinity  $E$ in $(ia_1,1,0,0)$ and contains $X$. But the line  $\bar E_1$ meets $X$ and has equation $-ia_1x+pw-ya_1^2=0, z+a_1^2w=0$ and hence meets $E$ where $w=z=0, x=ia_1y$,   which is  the same point.  Hence $a_1+it=0$ is the plane spanned by the two intersecting lines $X$ and $\bar E_1$.

It meets the surface $V$ in a third line: one that meets $X$ and $\bar E_1$. This can't be defined by a line through $x\in \PP^2$ since these are accounted for in the diagram so it must be derived from a conic through $x$. It cannot pass through $e_1$ since $x$ and $e_1$ get blown up and then  it would not meet $\bar E_1$ at all. It must  therefore be the line $C_1$, the proper transform of the conic missing $e_1$. So the divisor of $a_1+it=0$ is $X+\bar E_1+C_1$.

When $a_1+it=0$, $z-t^2=z+a_1^2$ which has, as we have seen, divisor $E+E_1+\bar E_1$.
Hence $z-t^2$ and $(a_1+it)(a_2+it)$ have a common divisor $\bar E_1+\bar E_2$. From the definition of $t$, $2X$ is also a common divisor. In homogeneous form we have the ratio of two cubic expressions:
$$f=\frac{(a_1(y-z+\sigma_2w)+i(x+i\sigma_1 z-i\sigma_3w))(a_2(y-z+\sigma_2w)+i(x+i\sigma_1 z-i\sigma_3w))w}{z(y-z+\sigma_2w)^2-w(x+i\sigma_1 z-i\sigma_3w)^2}.$$
When $w=0$ the denominator vanishes when $z=0$ or $(y-z)^2=0$ which is the divisor $E+2G$. The factor $w$ in the numerator  vanishes on $E+F+G$.  The divisor of the numerator can therefore be reduced to $ F+C_1+C_2$. 

As cohomology classes in $V$,
$(F+C_1+C_2)^2=-1-1-1+2+2=1$ and $-KC=1$ for any line $C$. By Riemann-Roch this means that $\dim H^0(V, {\mathcal O}(\bar D))=3$. A generic curve in the linear system has genus given by $2g-2 = K\bar D+\bar D^2 = -3+1$ and so $g=0$. Blowing up $(1,0,0,0)$ gives $D$ with $D^2=0$ and $\dim H^0(\hat V,{\mathcal O}(D)) = 2$, hence a pencil of rational curves.

Now describe the cohomology classes in terms of the blow-up of $\PP^2$. Let $H$ be the divisor of a line in $\PP^2$ then $H,E_1,\dots, E_4, X, F$ are generators of $H^2(V,\Z)$ and the divisor class $ D=F+C_1+C_2$ is given by 
$$D=4H-E_1-E_2+E_3+E_4-2X-F.$$ 
If $C_0$ denotes the line in $V$ defined by the conic in $\PP^2$ missing $f$ and $E_{34}$ the line joining $e_3$ to $e_4$ then we see that as divisor classes 
$$\bar D = F+C_1+C_2 \sim G+C_0+E_{34} \sim E+\bar E_3+\bar E_4.$$
These are three singular members of the pencil defined by  the same configuration of rational curves in $\hat V$. 
The last expression has the best interpretation in our context: it shows that  the linear system $E+\bar E_3+\bar E_4$ and its conjugate $E+ E_3+ E_4$  yield the two families of null curves defining the complexified conformal structure.
In fact the  denominator $z-t^2$ has divisor which is symmetric in the $a_i$ and is 
$$2X+E+2G+\bar E_1+\bar E_2+\bar E_3+\bar E_4$$
which reduces to $E+\bar E_3+\bar E_4$. This is the polar divisor of the  function $f$.
\end{prf}

\section{Compactifications}\label{comp}
The three concurrent lines $E,F,G$ on the cubic surface $V$  become  curves $E',F',G'$ of self-intersection $-2$ when we  blow up their common intersection to obtain $\hat V$. This gives a configuration of rational curves supporting an anticanonical divisor: $-K\sim 2C+E'+F'+G'$.

\hskip 4cm\begin{tikzpicture}
\draw[] (0,3) -- (6,3);
\draw[] (1,3) -- (1,0);
\draw[] (3,3) -- (3,0);
\draw[] (5,3) -- (5,0);
\node at (0.5,0) {${ E'}$};
\node at (2.5,0) {${ F'}$};
\node at (4.5,0) {${ G'}$};
\node at (4.5,0) {${ G'}$};
\node at (0.5,1.5) {$-{ 2}$};
\node at (2.5,1.5) {$-{ 2}$};
\node at (1.75,3.5) {${ C}$};
\node at (3.5,3.5) {${ -1}$};
\node at (4.5,1.5) {${ -2}$};
\end{tikzpicture}

In fact, each surface in the versal deformation of type $\Ek$ has a similar anticanonical configuration which fibrewise compactifies the simultaneous resolution:

\hskip 4cm\begin{tikzpicture}
\draw[] (0,3) -- (6,3);
\draw[] (1,3) -- (1,0);
\draw[] (3,3) -- (3,0);
\draw[] (5,3) -- (5,0);
\node at (0.5,0) {${ E'}$};
\node at (2.5,0) {${ F'}$};
\node at (4.5,0) {${ G'}$};
\node at (4.5,0) {${ G'}$};
\node at (0.25,1.5) {${ 3}-{ k}$};
\node at (2.5,1.5) {$-{ 2}$};
\node at (1.75,3.5) {${C}$};
\node at (3.5,3.5) {${- 1}$};
\node at (4.5,1.5) {${ -3}$};
\end{tikzpicture}

For the fibre over $Z_0$, $C$ comes from the line at infinity in $\PP^2/\Gamma$, resolving the three orbifold points of the action of $\Gamma$ on $\PP^1$, corresponding to stabilizers of vertices, edges and faces of the regular solid.  Moreover \cite{PB2},\cite{Pink} any nonsingular surface with such an anticanonical configuration is part of the simultaneous resolution of the corresponding singularity, so each fibre of $Z_+$ can be compactified in the same way to produce $\bar Z_+$. 
\begin{rmk} An ALE space has a conformal compactification as a self-dual orbifold and this configuration compactifies the twistor space as in \cite{PB2} and also Atiyah's paper on Green's functions \cite{At}.
\end{rmk} 
 A similar situation holds for $\Ak$ and $\Dk$ but we shall deal with these cases later. 

The $\C^*$-action on $\C^2/\Gamma$ is just the scalar action on $\C^2$ and an orbit in the resolution compactifies in $\bar Z_0$ to a projective line $L$ which meets a generic point on the exceptional curve $C$ transversally. This means $-KL = (2C+E'+F'+G') L=2$ and hence from the adjunction formula $L^2=0$ and the normal bundle of $L$ is trivial. 

Now $\bar Z_+\rightarrow \C$ is a fibration by compact surfaces so we can ask if $L\subset \bar Z_0$  deforms in the family. Because the normal bundle in $\bar Z_+$ is trivial, it does indeed deform into  a two-parameter family of such lines, giving a pencil in each fibre. At a point on $S_0$, of course fixed by  $\C^*$,  the action  on the tangent space is $(1,u,u)$, so a line in this family has an open set of points which lie in the  large Bialynicki-Birula stratum. Thus each deformation of $L$ is the compactification of a null curve  in the fibre $Z_u$ for  $u\ne 0$. 

The defining property of $L$ at $u=0$ is that it is transverse to $C$ in the compactification and  intersects $S_0$ transversally but intersects none of the other curves  in the resolution. This fixes its cohomology class. Then,

\begin{thm} Let $\bar Z_1$ be the compactification of the simultaneous resolution of a singularity of type $\Ek$ at $u=1$, $S_0$ the (real) central sphere and $C$ the (holomorphic) exceptional curve. 

Then the linear system whose  cohomology class  is dual via the intersection pairing to $[S_0]+[C]\in H^2(\hat Z_1,\Z)$, consists of a pencil of rational curves which intersect $Z_1\subset \bar Z_1$ in  null curves of the complexified conformal structure.
\end{thm}

\begin{rmk} In the  case of the moduli space of Higgs bundles, the analogue of the generic $\C^*$-orbit in $Z_0$ is the space of Higgs fields $H^0(\Sigma,\End V\otimes K)$ for a fixed stable bundle $V$. This is a Lagrangian submanifold and is the ``upward flow" of the $\C^*$-action from the point $[V]$ in the fixed point set.  There is a fibrewise compactification of the moduli space of $\lambda$-connections \cite{S1},\cite{S2}  and if $V$ is ``very stable", meaning it has no nilpotent Higgs field, then the vector space $H^0(\Sigma,\End V\otimes K)$ compactifies to a projective space. This deforms as we change the parameter $u$ and at $u=1$, in complex structure $J$, it becomes the projective completion of the affine space of moduli of flat  connections $\nabla_A$  such that the $(0,1)$-part $\bar\partial_A$ is equivalent to $V$. Together with its conjugate we obtain two transverse families of half-dimensional affine spaces which constitute the complexification of the K\"ahler metric on the moduli space of stable bundles. 

If there is a nilpotent Higgs field then the upward flow reaches another fixed point, like the $\C^*$-orbits whose limits are intersection points of the rational curves in the resolution. Note that, unlike the ALE case,  the twistor space for Higgs bundles cannot be compactified this way. 
\end{rmk}
\section{The divisor class}
 \subsection{Type ${\mathbf E}_{\mathbf k}$}
The compactification $\hat V$ in the $\DIV$ case has rational curves $E',F',G'$ of self intersection $-2$. Given that any anticanonical  configuration of rational curves as above defines the compactification in the versal family we can take as a model for $\Ek$ the $\DIV$ compactification blown up at  a real point $y$ on $G$  and further real points $e_5,\dots, e_{k-1}$ on $E$. We then achieve the required self-intersections $3-k,-2,-3$. This is a convenient substitute for the algebraic equation in $x,y,z$. 

Note that by blowing up real points we can extend the real structure and because those points are on the lines at infinity, from this point of view the sphere $S_0$ remains the same, as does its symplectic form. Only the conformal structure changes.

We observed in Proposition \ref{D4div} that the divisor $E+\bar E_3+\bar E_4$ defines the conformal structure for $\DIV$. We now have
\begin{thm} There is a choice of parameter $a\in \lie{h}$ such that the complexified conformal structure on the central sphere of type $\Ek$ is defined by a meromorphic function whose polar divisor is $E+\bar E_3+\cdots + \bar E_{k-1}$. 
\end{thm} 

\begin{rmk} The argument in Section \ref{comp} giving a pencil of rational curves is based on compactifying in $\PP^2/\Gamma$ the generic $\C^*$-orbit in $\C^2/\Gamma$. This gives a nonsingular curve of the pencil. The components $\bar E_5,\dots, \bar E_{k-1}$ correspond to the orbits whose limit is a fixed point on the line at infinity -- the $(k-3)$ vertices of the face of a regular solid. This is one  of the three singular members of the pencil. 
\end{rmk}
\begin{prf} Note that $(E+\bar E_3+\cdots + \bar E_{k-1})(E_3-E_2)=1$. The idea of proof  is to find a basis for the cohomology which is  orthogonal to $D=E+\bar E_3+\cdots + \bar E_{k-1}$ and $E,F,G$ and show that the $k-1$ generators together with $E_3-E_2$ as the trivalent vertex  intersect according to the relevant Dynkin diagram. 

Working as before with the blow-up of $\PP^2$ we write a general class as 
$$aH+\sum_{i=1}^{k-1}a_iE_i+bX+cY + dF.$$
Orthogonality to $F$ gives $d=0$ and to $E$ and $G$
\begin{equation}
a+b+c=0=2a+\sum_{i=1}^{k-1} a_i.
\label{ortho1}
\end{equation}
The intersection with $D=E+\bar E_3+\cdots + \bar E_{k-1}$ is 
\begin{equation}
(k-3)(a+b)+\sum_{i=3}^{k-1}a_i.
\label{ortho2}
\end{equation}
Clearly $E_2-E_1$ satisfies  (\ref{ortho1}) and from (\ref{ortho2}) is orthogonal to $D$, similarly    $E_4-E_3,\dots, E_{k-1}-E_{k-2}$. Together with $E_3-E_2$ this gives a chain of classes of self-intersection $-2$.

For $b=-1$ and $a=1$ we have $P=H-E_1-E_2-X$, and then $P^2=-2$ and $P(E_3-E_2)=-1$. 

\noindent 1. When $k=6$, $b=-1$ and $a=2$ yields $Q = 2H-E_2-E_3-E_4-E_5-X-Y$. Then  $Q^2=-2$ and 

\hskip 4cm \ddei Q.{-P}.E_2-E_1.E_3-E_2.E_4-E_3.E_5-E_4.

gives the $\EVI$ Dynkin diagram.

\noindent 2. When $k=7$,  $b=-1$ and $a=2$ gives  $R=2H-(E_3+E_4+E_5+E_6)-X-Y$ with $R^2=-2$ and $PR=1$. Then 

\hskip 4cm \ddeii P+R.{R}.E_2-E_1.E_3-E_2.E_4-E_3.E_5-E_4.E_6-E_5.

gives the $\EVII$ Dynkin diagram.

\noindent 3. When $k=8$,  $b=-2$ and $a=3$ gives $S=3H-(E_2+E_3+E_4+E_5+E_6+E_7)-2X-Y$ with $S^2=-2$, $S(E_3-E_2)=0$, $S(E_2-E_1)=1$ and $SP=0$,  from which we have the $\EVIII$ Dynkin diagram:

\hskip 3cm \ddeiii S.{-P}.E_2-E_1.E_3-E_2.E_4-E_3.E_5-E_4.E_6-E_5.E_7-E_6.
\end{prf}
\subsection{Periods} 
From Kronheimer's approach, the moduli are given by the periods of the holomorphic $2$-form (only the real part has non-zero periods). We have a basis above associated to simple roots of the $\Ek$ root system, which places the moduli in the Cartan subalgebra $\lie{h}$.  Each generator is a sum of terms of the form $A-B$ where $A,B$ are represented by exceptional curves. For $E_{i+1}-E_i$ this is obvious but, recalling that $E_{ij}$ denotes the  proper transform of the line joining $e_i$ and $e_j$, we have 
$$P=E_{12}-X,\qquad Q=(E_{23}-X)+(E_{45}-Y)$$
$$ R=(E_{34}-X)+(E_{56}-Y),\qquad S=(E_{23}-X)+(E_{45}-X)+(E_{67}-Y). $$
The curves $E_i$ intersect $E$ but not $G$. A curve $E_{ij}$ does not intersect $E$ since the line joining $e_i$ to $e_j$ in $\PP^2$ is not tangential to the conic  defining $E$. Also $X,Y$ intersect $G$ but not $E$. Thus each generator is of the form $A-B$ where $A,B$ intersect $E$ or $G$, but not both. 

The holomorphic symplectic form on $Z_1$ has a simple pole on $E$ and $G$ and we can calculate its period on $A-B$ from the Poincar\'e residue, a one-form on the polar divisor. If $A$, $B$ intersect $E$, say, in $a$ and $b$ then take a path $\gamma$ from $a$ to $b$ in $E$ and a circle bundle in the normal bundle of $E$ over $\gamma$. Then $A-B$ is homologous to the cycle obtained by cutting out discs around $a,b$ in $A,B$ and connecting with this cylinder. Since a holomorphic 2-form vanishes on a complex curve, the period on $A-B$ is the  integral over the cylinder which is $2\pi i$ multiplied by the integral of the Poincar\'e residue along $\gamma$.  

In our situation, $Z_1$ is obtained by blowing up points on $\PP^2$ and the holomorphic form is the transform of  a 2-form on $\PP^2$ with a pole on a singular cubic curve consisting of a conic (defining $E$) and a tangent to it (defining $G$). By a projective transformation we can take the conic to be $yz=x^2$ and the line to be the line at infinity $z=0$. Then the holomorphic form in affine coordinates $x,y$ on  $\C^2$ is a (possibly complex) multiple of 
$$\omega=\frac{1}{(y-x^2)}dx\wedge dy=\frac{1}{(y-x^2)}dx\wedge d(y-x^2)$$
so the residue on $E$ is $dx$ on the parabola $y=x^2$. Thus $x$ is a natural parameter on $E$.

As for $G$, near the line at infinity using coordinates $\tilde y=y/x, \tilde z =z/x$ we have 
$$\omega=-\frac{1}{\tilde y(\tilde y \tilde z-1)}d\tilde z\wedge d\tilde y$$
giving the residue $d(y/x)$ on $Z=0$ and $\tilde y=y/x$ the natural parameter there: the line $y=mx+c$ in $\C^2$ intersects the line at infinity at $\tilde y=m$. The point $f=E\cap G$ is $(0,1,0)$, $\tilde y=\infty$.  We can choose the point $x\in G$ to be $\tilde y=0$.

From this point of view  $e_i\in E$ is defined by a point $(x,y)=(a_i,a_i^2) \in \C^2$ and so the period of $E_{i+1}-E_i$ is $2\pi i(a_{i+1}-a_i)$. The line $E_{ij}$ joining $e_i$ to $e_j$ has slope $(a_j^2-a_i^2)/(a_j-a_i)=a_j+a_i$, so if  $y\in G$ be given by $\tilde y= a_0$, we have the following periods for $P,Q,R,S$:
$$\omega(P)=2\pi i(a_1+a_2),\qquad \omega(Q)=2\pi i(a_2+a_3+a_4+a_5-a_0)$$
$$\omega(R)=2\pi i(a_3+a_4+a_5+a_6-a_0),\quad \omega(S)=2\pi i(a_2+a_3+a_4+a_5+a_6+a_7-a_0).$$
From these formulae one can derive the positivity conditions for the vertices of the Dynkin diagram to represent holomorphic curves -- the necessary condition for our ALE metric.

\begin{rmk} We can replace the datum of the point $y\in G$ by the other tangent from $y$ to the conic which meets it at a point $(a_k,a_k^2)$. Then the tangent direction is $a_0=2a_k$ and so the parameters consist of $k$ points on $E$. The birational map on $\PP^2$ given by $(x_0,x_1,x_2)\mapsto (x_0^2,x_0x_1,x_1x_2)$ sends the conic $(1,t,t^2)$ to the cuspidal cubic $(1,t,t^3)$. This makes our approach to the $\Ek$ simultaneous resolution compatible with the traditional approach of $k$ points on the cubic \cite{Ty}. Geometrically, $G$ becomes an exceptional curve when $y\in G$ is blown up and blowing it down, the conic $E$ becomes the singular cubic.
\end{rmk}
\subsection{Type $\Dk$}
The compactification for $\Dk$ is the following

\hskip 4cm\begin{tikzpicture}
\draw[] (0,3) -- (6,3);
\draw[] (1,3) -- (1,0);
\draw[] (3,3) -- (3,0);
\draw[] (5,3) -- (5,0);
\node at (0.5,0) {${ E'}$};
\node at (2.5,0) {${ F'}$};
\node at (4.5,0) {${ G'}$};
\node at (4.5,0) {${ G'}$};
\node at (0.25,1.5) {$-{ k-2}$};
\node at (2.5,1.5) {${ -2}$};
\node at (1.75,3.5) {${C}$};
\node at (3.5,3.5) {${ -1}$};
\node at (4.5,1.5) {${ -2}$};
\end{tikzpicture}

and can clearly be achieved by blowing up further points on $E$. The divisor has the same form as the $\Ek$ case:
$$D=E+\sum_{i=3}^k\bar E_i.$$
\subsection{Type $\Ak$}
For the sake of completeness we finally approach the $\Ak$ case using the same method. 
As before we deal with $k=2\ell-1$. Then the compactification is 

\hskip 5cm\begin{tikzpicture}
\draw[] (0,3) -- (4,3);
\draw[] (1,3) -- (1,0);
\draw[] (3,3) -- (3,0);
\node at (0.5,3.5) {${C}$};
\node at (0.5,1.5) {${- \ell}$};
\node at (2.7,1.5) {${- \ell}$};
\node at (2,3.4) {${ 0}$};
\node at (0.5,0) {${ D_1}$};
\node at (2.6,0) {${ D_2}$};
\end{tikzpicture}

A model for the surface  $\bar Z_1$ is to take the projective bundle $\PP({\mathcal O}(\ell)\oplus 1)$ with the infinity section $D_1$ and blow up $2\ell$ points on the zero section $D_0$ with $D_0^2=\ell$ to obtain $D_2$. A general fibre is then $C$ and as divisor classes $D_2=\ell C+D_1-\sum_i E_i$ where $E_i$, $1\le i\le 2\ell$ are the exceptional curves. The anticanonical class is $-K=2C+D_1+D_2$ from the compactification picture but is also 
$$-K=(\ell+2)C+2D_1-\sum_{i=1}^{2m}E_i$$
from the model. 

Now $F_i=C-E_i$ satisfies $F_i^2=-1$, $-K.F_i=1$ and so $H^0(\bar Z_1, {\mathcal O}(F_i))=1$ and defines another exceptional curve since $C=E_i+F_i$. The projection onto $\PP^1$ then has singular fibres a pair of intersecting exceptional curves $E_i,F_i$ over points $a_1,\dots, a_{2\ell}\in \PP^1$. 

To determine the divisor class, we note that the orthogonal complement  of $C,D_1$ and $D_2$ is spanned by cohomology classes $E_2-E_1, E_3-E_2,\dots, E_{2\ell}-E_{2\ell-1}$ which intersect in a chain giving the $\Ak$ Dynkin diagram and the  class which picks out the central element is 
$$D=D_1+\sum_{i=1}^{\ell}F_i.$$
To relate this to the explicit meromorphic function in Section \ref{Acentral} observe that $\ell C+D_1+D_2-\sum_iE_i\sim 2D_2$ and since $D_1^2=-\ell$ there is, up to a multiple,  a unique section $s_1$ of ${\mathcal O}(D_1)$ and hence a unique section $x$ of $\ell C+D_1+D_2$ vanishing on the exceptional curves $E_i$. Similarly a section $y$ vanishing on $F_i$. Let $z$ be an affine coordinate on $\PP^1$ and $z-a_i$ denote the section of ${\mathcal O}(1)$ whose fibre is $E_i+F_i$ then we have the equation of the affine surface 
$$xy=\prod_{i=1}^{2\ell}(z-a_i).$$
Then $x$ and $(z-a_1)\dots(z-a_{\ell})$ have common divisors $E_1,\dots, E_{\ell}$ and $D_2$, which reduce the zeros of $x/(z-a_1)\dots(z-a_{\ell})$ to $\sum_{\ell+1}^{2\ell}E_i+D_2$ in the numerator and $D_1+\sum_1^{\ell} F_i$ in the denominator as required.
\section{ALF deformations}
The multi-instanton hyperk\"ahler metrics have deformations called multi-Taub-NUT spaces. These are obtained from the Gibbons-Hawking Ansatz by adding a  constant $1$ to the potential function $V$. They are no longer asymptotically locally Euclidean but instead have a different decay property at infinity described as being asymptotically  locally flat (ALF). 

Kronheimer \cite{PB0} observed that one could interpret this space as the moduli space of charge $1$ $SU(2)$-monopoles on $\R^3$ in a field of $k+1$ Dirac monopoles located at points ${\mathbf a_i}\in \R^3$. This modification also holds in the $\Dk$ case (though not for $\Ek$) and the paper \cite{CH} was motivated by the interpretation as a moduli space of charge $2$ monopoles in a field of Dirac monopoles, generalizing in some way the so-called Atiyah-Hitchin metric \cite{AH} where there are no  Dirac singularities. 

The twistor spaces for the ALE and ALF versions have isomorphic pieces $Z_+, Z_-$ but patching them together involves a certain exponential function which means that the whole twistor space cannot be compactified. Nevertheless, when the parameters are chosen so that the  metric admits a circular symmetry we may consider how the metric structure of the central sphere changes. As we have seen, the symplectic structure is determined just by considerations of $Z_+$, so it is the conformal structure which changes. 

Consider first the $\Ak$ case. The passage from ALE to ALF is to change $V$ to $1+V$. The $(1,0)$ forms of the conformal structure then change from $Vdz+id\theta$ to 
$$(1+V)dz+id\theta=Vdz+id(\theta-iz).$$
The symplectic form is $d\theta\wedge dz$ so this is the transform by a complex symplectic transformation  $(z,\theta)\mapsto (z, \theta -iz)$ which is obtained by integrating the Hamiltonian vector field for the function $iz^2/2$. 

In terms of the affine surface $xy=\prod_i(z-a_i)$ this is the transformation 
$(x,y,z)\mapsto (e^zx,e^{-z}y,z).$ It is a well-defined holomorphic but non-algebraic symplectic transformation of the surface and transforms the complex null curves of the ALE 
space to the ALF version. The conformal structure is obtained by taking this family together with its conjugate. Note that the real structure here is $(x,y)\mapsto (\bar y,\bar x)$ so this is not a real transformation, indeed that would give a diffeomorphic metric.  

This exponential expression is reflected in the patching together of $Z_+$ and $Z_-$ to form the twistor space. The analogous procedure as described in \cite{CH} shows that the $\Dk$ surface 
$$x^2-zy^2=-\frac{1}{z}\left(\prod_{i=1}^k(z+a_i^2)-\prod_{i=1}^ka_i^2\right)+2y\prod_{i=1}^k a_i$$
is transformed to ALF form by integrating the Hamiltonian vector field for the function $iz$.

\vskip 1cm
 Mathematical Institute,  Woodstock Road, Oxford OX2 6GG, UK
 
 hitchin@maths.ox.ac.uk

 \end{document}